\documentclass[a4paper,11pt,twoside]{article}	
\usepackage{amsmath,amssymb,amsthm,amsfonts,color,graphicx}
\usepackage{subcaption}
\usepackage{setspace}
\usepackage[margin=3cm]{geometry}
\newtheorem{theorem}{Theorem}
\newtheorem*{theorem*}{Theorem}
\DeclareMathOperator{\hr}{\mathbb H^2\times\mathbb R}

\DeclareMathOperator{\hh}{\mathbb H^2}
\DeclareMathOperator{\rr}{\mathbb R}

\newcommand{\N}{\mathbb{N}}

\newcommand{\R}{\mathbb{R}}

\newcommand{\s}{\mathbb{S}}
\newcommand{\h}{\mathbb{H}}

\newcommand{\B}{\mathbb{B}}

\newtheorem{proposition}{Proposition}

\theoremstyle{definition}
\newtheorem{claim}{Claim}

\numberwithin{equation}{section}

\begin{document}

\begin{title}
  {Slab Theorem and Halfspace Theorem for constant mean curvature
    surfaces in $\hr$}
\end{title}
\vskip .2in

\begin{author} {Laurent Hauswirth, Ana Menezes and Magdalena Rodr\'\i
    guez\thanks{Research partially supported by MCIN/AEI/10.13039/501100011033/ grant no. PID2020-117868GB-I00, FEDER/Andalucía grants no.\
A-FQM-139-UGR18 and P18-FR-4049 and by MINECO/FEDER grant no.\
MTM2017-89677-P.}}

\end{author}

\newcommand{\Addresses}{{
  \bigskip
  \footnotesize

\textsc{D\'epartement de Math\'ematiques
  Universit\'e de Marne-la-Vall\'ee, Cit\'e Descartes, France}\par\nopagebreak
  \textit{E-mail address:} \texttt{hauswirth@math.univ-mlv.fr}

  \medskip

   \textsc{Department of Mathematics, Princeton University, USA}\par\nopagebreak
  \textit{E-mail address}: \texttt{amenezes@math.princeton.edu}

  \medskip

 \textsc{Departamento de Geometr\'\i a y Topolog\'\i a, Universidad
  de Granada, Spain}\par\nopagebreak
  \textit{E-mail address}: \texttt{magdarp@ugr.es}

}}

\date{}

 \maketitle

\begin{abstract}
We prove that a properly embedded annular end of a surface in $\hr$
with constant mean curvature $0<H\leq \frac{1}{2}$ can not be
contained in any horizontal slab. Moreover, we show that a properly
embedded surface with constant mean curvature $0<H\leq \frac{1}{2}$
contained in $\hh\times[0,+\infty)$ and with finite topology is
  necessarily a graph over a simply connected domain of $\hh$.  For
  the case $H=\frac{1}{2}$, the graph is entire.

\end{abstract}

\begin{flushleft}
\end{flushleft}

\begin{flushleft}
\end{flushleft}

\section{Introduction}

The theory of constant mean curvature (CMC) $H>0$ surfaces in $\hr$
drew a lot of attention after the work by Abresch and
Rosenberg~\cite{AR}, where they described a Hopf-type holomorphic
quadratic differential on any such surface, and characterized the CMC
rotational spheres for $H>\frac 12$ as the only immersed CMC spheres
in this space~\cite{AR, HH, NR06,NSST}.  For $0\leq H\leq \frac 12$,
there are no compact CMC examples. This is why $H=\frac{1}{2}$ is
called the critical value for the mean curvature in $\hr$. The CMC
rotational simply connected examples for $0< H\leq \frac 12$ are
entire graphs of paraboloid-type shape (see Section \ref{sec-parab}
for more details).  The geometric behavior of CMC surfaces in $\hr$
for $H>\frac 12$ is, in some sense, analogous to the one of surfaces
of positive CMC in $\R^3$. For instance, for these values of the mean
curvature there exist spheres and a 1-parameter family of annuli
invariant by a vertical translation similar to the Delaunay's examples
(see \cite{NSST} and references therein).

An important theorem by Hoffman and Meeks~\cite{HM} in the classical
theory of minimal surfaces in $\R^3$ is the Halfspace Theorem saying
that there are no properly immersed non-flat minimal surfaces in a
halfspace. However, in $\hr$ this result does not hold: there are many
entire minimal graphs and rotational annuli (called catenoids)
contained in a slab constructed by Nelli and Rosenberg~\cite{NR02}.
There are also other properly embedded minimal annuli of bounded
height constructed in~\cite{FMMR}.  Notice that the existence of
spheres for $H>\frac 12$ and paraboloid-type graphs for $0< H\leq
\frac 12$ implies that entire graphs of bounded height cannot exist in
$\hr$ for $H>0$ by the maximum principle.

In \cite{CHR} Collin, Hauswirth and Rosenberg proved that a properly
embedded simply connected minimal surface in a slab $S$ of height less
than $\pi$ must be an entire graph. More generally, they proved that
each end of a minimal surface properly embedded in $S$ with finite
topology is a graph outside a compact domain.  For the case of CMC
surfaces with $0<H\leq \frac 1 2$ the behavior is different. We prove
that there are no examples in a slab. 

\begin{theorem*}[Slab Theorem]
  Let $M\subset \hr$ be a surface (possibly with boundary) with
  constant mean curvature $0<H\leq \frac{1}{2}$ and at least one
  properly embedded annular end. Then $M$ can not be contained in a
  horizontal slab $\h^2\times[0,L]$, for any $L>0$.

In particular, there are no properly embedded CMC surfaces for $0<H\leq
\frac{1}{2}$ with finite topology contained in a horizontal slab of
$\hr$.
\end{theorem*}

For any $0<H\leq \frac 12$, Manzano and Torralbo constructed
in~\cite{MT1} properly immersed CMC surfaces contained in a slab.
These examples are invariant by a group of symmetries induced by a
tessellation of $\mathbb{H}^2$ by regular polygons. A fundamental
domain of any of these examples is compact and its lift in $\hr$
has, by periodicity, only one end of infinite topology.  If these
examples were embedded, it would show that the hypothesis of having an
annular end in the theorem above is sharp.

For CMC $H=\frac 1 2$, there exists a halfspace theorem~\cite{HRS} for
complete embedded CMC surfaces in $\hr$ lying on one side of a
horocylinder (with some condition on the mean curvature vector); the
only such examples are parallel horocylinders. Nelli and Sa
Earp~\cite{NS} also proved that the only CMC surfaces with $H=\frac
12$ contained in the mean convex side of the rotationally invariant
paraboloid-type entire graph are translated copies of the graph.  If
we think of CMC surfaces on one side of a horizontal slice, the only
known examples are entire graphs. Here we prove that, for $0<H\leq
\frac{1}{2}$, the only properly embedded CMC surfaces with finite
topology contained in one side of a horizontal slice are graphs.

\begin{theorem*}[Halfspace Theorem]
  Let $M\subset\hh\times[0,+\infty)$ be a properly embedded surface
    with constant mean curvature $0<H\leq \frac{1}{2}$ and finite
    topology. Then $M$ is necessarily a graph over a simply connected
    domain of $\hh$.  For $H=\frac{1}{2}$ the
    graph is entire.
\end{theorem*}

\section{Preliminaries}

In this section we will set up some notations and introduce some
classes of constant mean curvature (CMC) graphs in $\hr$ that we will
use as barriers.  Throughout this paper we consider the cylinder model
for $\hr$. We consider $\mathbb H^2=\{(x,y)\in\rr^2;\ x^2+y^2<1\}$
endowed with the hyperbolic metric
$g_{-1}=\frac{4}{(1-x^2-y^2)^2}g_0,$ where $g_0$ denotes the Euclidean
metric in $\rr^2.$ We will then consider the unit solid cylinder with
the product metric $g=g_{-1}+dt^2$ as model for $\hr$. In this model
there is a natural notion of asymptotic boundary of $\hr$ where
$(\partial_\infty\h^2)\times\R=\s^1\times\R$.

\subsection{Constant mean curvature Scherk graphs}
\label{scherk}

For any $H\in(0, \frac{1}{2}],$ Hauswirth, Rosenberg and Spruck
  \cite{HRS2} described necessary and sufficient conditions over a
  compact admissible domain in order to guarantee the existence of a
  graph with constant mean curvature $H$ assuming infinite boundary values. A compact domain
  $\Omega$ in $\hh$ is said to be admissible if it is simply connected
  and its boundary $\partial\Omega$ is a polygon with sides $\{A_i\}$
  and $\{B_i\}$, with neither two consecutive $A_i$ edges nor $B_i$
  edges, and all satisfying $\kappa(A_i)=2H$ and $\kappa(B_i)=-2H$,
  where $\kappa$ denotes the geodesic curvature with respect to the
  interior of $\Omega$. The Jenkins-Serrin problem they considered
  consists of finding a solution to the equation for CMC graphs in
  $\Omega$ of mean curvature $H$ (we will call it $H$-graph) which
  assumes boundary values $+\infty$ on each $A_i$ and $-\infty$ on
  each $B_i$. In the case the domain is an ideal quadrilateral with
  two opposite $A_1, A_2$ edges and two opposite $B_1, B_2$ edges,
  Folha and Melo constructed complete examples $-$ which are called
  complete Scherk $H$-graphs $-$ for any $0<H< \frac{\sqrt{2}}{2}$
  (see Appendix in \cite{FoM}). Since we will use as barriers these
  complete Scherk $H$-graphs for $0<H<\frac{1}{2}$, here we show their
  existence for this wider range of the mean
  curvature.

\begin{proposition}
For any $0<H<\frac{1}{2}$, there exists a complete Scherk $H$-graph
over an ideal domain, called Scherk domain, bounded by curves $A_i$
and $B_i$ (as above) with $\kappa(A_i)=2H$ and $\kappa(B_i)=-2H,$ for
$i=1,2.$
\label{scherk}
\end{proposition}

\begin{proof}
We use a Plateau conjugate method.  Melo \cite{Melo} constructed a
complete minimal Scherk-type graph over an ideal quadrilateral in any
$\widetilde{PSL}(2,\mathbb R)$ space. Using Daniel's correspondence
\cite{D} and the techniques by Castro-Infantes, Manzano and the third
author in~\cite{CMR}, these graphs correspond to complete Scherk
$H$-graphs into $\hr$ as desired, with ${0<H<\frac 1 2}$.
\end{proof}

\subsection{Bigraph horizontal annuli with constant mean curvature $H=\frac 1 2$}

\label{sec.annuli}

\begin{figure}[h]
   \centering \includegraphics[height=4cm]{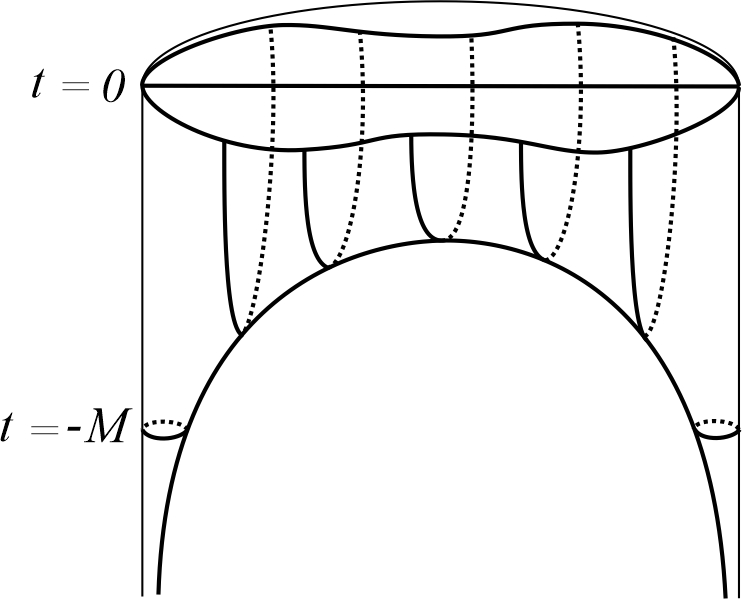}
\caption{Half of an horizontal catenoid $\mathcal{C}{_a}$ which is a
  graph ${\cal G}_a$ over a domain in
  $\hh\times\{0\}$.}
\label{figure20}
\end{figure}

For $H=\frac 1 2$, instead of using Scherk-type graphs as barriers we
will use the one-parameter family of horizontal annuli
$\{\mathcal{C}{_a}\}_{a>0}$ constructed by Daniel and the first author
in~\cite[Section 8]{DH}, called horizontal catenoids, whose boundary
at infinity consists of two vertical lines (see also~\cite{CMR} for an
alternative construction).  Up to an isometry we can assume that any
horizontal catenoid ${\cal C}_a$ is symmetric with respect to the
vertical plane $\{y=0\}$ that separates both ends and with respect to
the horizontal plane $\hh\times\{0\}$, and the lower half of ${\cal
  C}_a$ is a graph which we denote by ${\cal G}_a$ (see Figure
\ref{figure20}).  Each graph $\mathcal{G}_a$ is isometric to a minimal
surface $\Sigma_a$ in the Heisenberg space Nil$_3$, called its
conjugate surface, and share the same values of the angle function
$\nu=<N,\partial_t>$ (see Daniel's correspondence in~\cite{D}). This
minimal surface $\Sigma_a$ is a helicoidal surface bounded by two
vertical geodesics. These boundary geodesics correspond to the
horizontal curves at height~$0$ of ${\cal G}_a$, where $\nu =0$.

The parameter of the family $\{\mathcal{C}{_a}\}_{a>0}$ corresponds to
the size of the neck of the annuli, where by {\it neck} we mean the
(compact) intersection curve between the annulus and the vertical
plane of symmetry $\{y=0\}$.

When $a$ goes to zero, the limit surface consists of the union of two
tangent horocylinders (a pinching is produced in this case). The limit
surface become vertical everywhere, even if we consider different
translated copies of the annuli, so the limit domain of the domains
where the graphs ${\cal G}_a$ are defined is the union of two tangent
horodisks and is folliated by divergence lines (sets of points where
the gradient of the functions are unbounded), that are horocycles at
the same two points at infinity.

\begin{figure}[h]
   \centering \includegraphics[height=3.5cm]{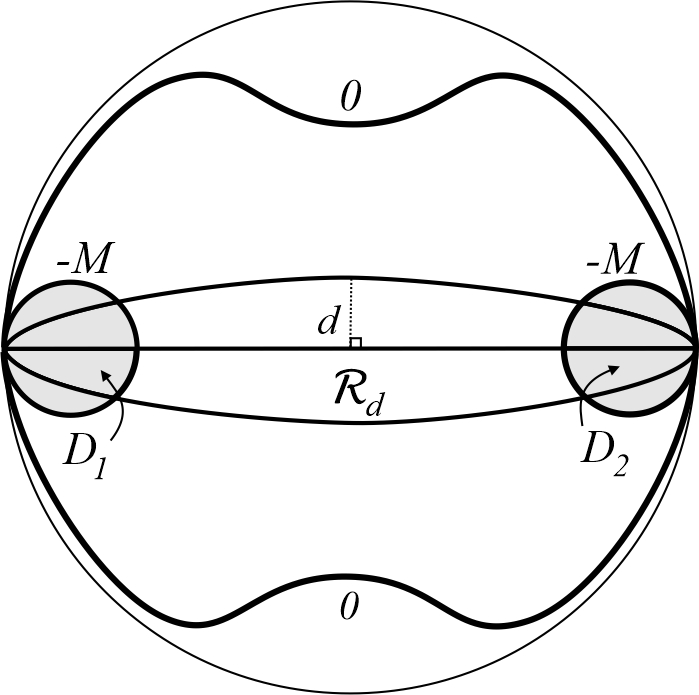}
   \caption{Projection of the graph ${\cal G}_a$ and the domains $D_1$
     and $D_2$ over which $\mathcal {G}_a\cap\{t<-M\}$ projects, for
     large $M$.  }
\label{figure21better}
 \end{figure}

Now translate the catenoids so that any ${\cal C}_a$ is tangent to
$\hh\times\{0\}$ at the origin of $\hh$ and this point is contained in
the neck of the annulus.  When $a$ diverges, the necks of the annuli
(all of them passing through the origin) become as large as we want.
Hence, when $a$ goes to $+\infty$, the graphs ${\cal G}_a$ converge to
the entire graph ${\cal I}$ given explicitly by Sa Earp~\cite[equation
  (31)]{Sa} which is invariant by a one-parameter family of hyperbolic
translations. Notice that in~\cite[Remark 3.7]{CMR} it is proved that
the conjugated minimal surfaces of ${\cal G}_a$ in Nil$_3$ converge to
the minimal entire graph in Nil$_3$ invariant by the isometric
translations along a horizontal geodesic, and Daniel proved
in~\cite[Example 5.6]{D} that the conjugate surface of this entire
minimal graph is the graph ${\cal I}$. In particular, for any $d>0$,
there exists $a_0$ such that ${\cal R}_d$ is contained in the domain
where the graph ${\cal G}_a$ is defined for any $a\geq a_0$, where
${\cal R}_d\subset\hh$ denotes the region bounded by the two
equidistant curves at distance $d$ to the horizontal geodesic
$\{y=0\}\subset\hh$, see Figure \ref{figure21better}.

On the other hand, the conjugate surface $\Sigma_a\subset$Nil$_3$ is
foliated by straight lines with a Gauss map which is horizontal at
infinity. Hence, on each straight line the function $\nu$ converges
uniformly to zero at infinity. This shows that, outside a compact set,
the unit normal vector at any point of the graph $\mathcal{G}_{a}$ is
arbitrarily close to horizontal. In particular, using different slide
back sequences we can prove that, for $M>0$ large enough, the
translation of the surface $\mathcal{G}_{a}\cap \{t< -M\}$ converges
to a horocylinder. This uniform convergence proves that the horizontal
curves $\mathcal{G}_{a}\cap \{t=-M\}$ are close to two horocycles and
contain any half equidistant curves to a geodesic having the same
points at infinity (corresponding to the ends of the horizontal
catenoid). Up to an isommetry we can assume that this geodesic is
$\{y=0\}\subset\hh$. In particular, $\mathcal {G}_a\cap\{t<-M\}$ is a
graph over two unbounded domains $D_1$ and $D_2$ so that ${\cal
  R}_d\setminus(D_1\cup D_2)$ is compact for any $d>0$.


\subsection{Entire constant mean curvature graphs in $\hr$}
\label{sec-parab}

There is a well known class of {\it entire} CMC graphs for any $H\in
(0, \frac{1}{2}]$ that are rotationally invariant complete vertical
  graphs with empty asymptotic boundary in
  $(\partial_\infty\h^2)\times\R$ (see, for
  instance,~\cite{NR06,NSST}). We are going to call them {\it
    paraboloids}. One such surface is the graph of a convex function $u$ that
  diverges to $+\infty$ when approaching $\partial_\infty\h^2$ and
  takes its global minimum at height $0$ (up to a vertical
  translation). We will denote by $\cal P^+$ this surface throughout
  this paper. The symmetric surface
  with respect to the horizontal slice at height $0$ will be denoted
  by ${\cal P}^-$.

  When $H= \frac{1}{2}$, ${\cal P}^+$ is the graph of the function (in
  polar coordinates)
$$
u(r,\theta)=\frac{1}{\sqrt{1-r^2}},\quad \mbox{where }\ 0\leq r<1.
$$

\begin{figure}[h]
   \centering \includegraphics[height=2.85cm]{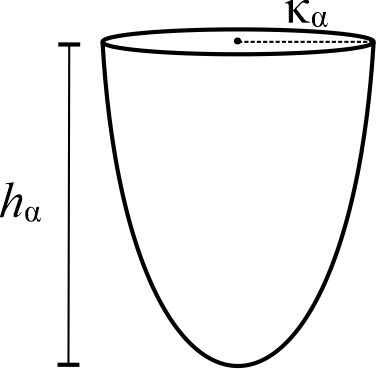}
   \caption{ ${\cal P}^+\cap\left( \hh\times\left[0, h_\alpha\right]\right)$}
  \label{paraboloid}
 \end{figure}

  The mean curvature vector of the paraboloid $\mathcal P^+$ points
  upwards, and we orient the surface by the unit vector field
  $N^+$ such that $\nu^+=\langle N^+,\partial_t\rangle$ is
  positive. Since ${\cal P}^+$ is rotationally invariant, we can think
  of $\nu^+$ as a function on $r\in[0,1)$ in polar coordinates. Moreover, since
    ${\cal P}^+$ is the graph of a convex function whose tangent plane
    is becoming vertical at infinity, $\nu^+$ is strictly decreasing and
    takes all values in $(0,1]$.  In particular, for any
  $\alpha\in(0,1)$ there exist $k_\alpha,h_\alpha>0$ such that the
  region of $\mathcal{P}^+$ where $\nu^+\geq \alpha$ coincides with
  ${\cal P}^+\cap\left( \hh\times\left[0, h_\alpha\right]\right)$ and is bounded by a horizontal circle
  of radius $k_\alpha$ (see Figure \ref{paraboloid}). We observe that both $h_\alpha$ and $k_\alpha$
  diverge as $\alpha$ goes to $0$.

The paraboloid ${\cal P}^+$ (as well as ${\cal P}^-$) separates the
ambient space into two connected components. We denote by ${\cal
  P}_{int}$ the mean-convex component (the one where the mean
curvature vector points to) and by ${\cal P}_{ext}$ the other one. We
will use this notation in the proof of the Slab Theorem.

There are many other examples of complete CMC graphs in $\hr$.  In
fact, for any $H\in (0, \frac{1}{2}]$ there are families of CMC graphs
  invariant by a hyperbolic or a parabolic translation, not
  necessarily entire when $H<\frac 12$ (see~\cite{Sa} or the appendix
  in~\cite{MRR11}). The intersection of the asymptotic boundary of all
  these examples with $(\partial_\infty\h^2)\times\R$ is contained in
  one or two vertical lines, and some of these examples are contained
  in one side of a horizontal slice.

Moreover, there exist many CMC $H=\frac 1 2$ entire graphs which are
not rotationally invariant, obtained as deformation of the paraboloid,
with empty asymptotic boundary in $(\partial_\infty\h^2)\times\R$ (see
\cite{CH} for more details).

\section{Slab Theorem}
\label{sec:slab}

Collin, Hauswirth and Rosenberg~\cite{CHR} have proved that an annular
end of a properly immersed minimal surface contained in a slab of
height less than $\pi$ is a multigraph outside a compact domain with a
finite number of sheets. When $0<H\leq \frac 12$, the following
theorem proves the non-existence of a properly embedded annular end
with CMC $H$ in a slab of any height $L$. We will denote by ${\cal
  S}_L$ the horizontal slab $\h^2\times[0,L]$, for any $L>0$.

\begin{theorem}
  Let $M\subset \hr$ be a surface (possibly with boundary) with
  constant mean curvature $0<H\leq \frac{1}{2}$ and at least one
  properly embedded annular end. Then $M$ can not be contained in
  ${\cal S}_L$, for any $L>0$.

In particular, there are no properly embedded CMC surfaces for
$0<H\leq \frac{1}{2}$ with finite topology contained in a horizontal
slab of $\hr$.
\label{thm1}
\end{theorem}

\begin{proof}
Let us suppose, by contradiction, that there exists one such surface
$M$ contained in ${\cal S}_L$ and call $E$ a properly embedded annular
end of $M$. Then $E\subset{\cal S}_L$ is an annulus with compact
boundary~$\partial E$. Since $E$ is properly embedded, there exists a
compact disk $D$ (not necessarily minimal) with $\partial D=\partial
E$ such that $E\cup D$ is a surface that separates $\hr$ into two
connected components. Along $E$, the mean curvature vector
$\overrightarrow{H}$ distinguishes these two components. We call
interior component the one where $\overrightarrow{H}$ points to and
exterior component the other one.  We denote by $N$ the unit normal
vector to $E$ such that $\overrightarrow{H}= H N$, and $\nu=\langle
N,\partial_t\rangle$ denotes what we call the angle function of $E$.

Let us consider a paraboloid $\mathcal P^+$ with the same mean
curvature $0<H\leq \frac 12$ as $E$ (see Section \ref{sec-parab}).
Its mean curvature vector is pointing upwards, and we orient the
paraboloid by the unit vector field $N^+$ such that $\nu^+=\langle
N^+,\partial_t\rangle$ is positive.  Since $\nu^+$ takes all values in
$(0,1]$ and ${\cal P}^+$ is rotationally invariant, for any point
  $p\in E$ with $\nu(p)>0$, we will be able to translate ${\cal P}^+$
  in such a way that it passes through $p$ with $N^+(p)=N(p)$. In the
  case $\nu(p)<0$, we instead consider a translation of ${\cal P}^-$
  to find a paraboloid with the same mean curvature vector as $E$ at
  $p$.


  
    \begin{claim}\label{cl}
    There exists $\alpha_0 \in (0,1)$ such that $\{p \in
    E\ :\ |\nu(p)| \geq \alpha_0\}$ contains a sequence of diverging
    points $(p_n)_{n \in \N}$, i.e.,
     \begin{equation}
      d(p_n):= \mbox{dist}_{\h^2}\left(\pi(p_n),\pi(\partial E)\right)
       \to +\infty \hbox{ and } |\nu (p_n)| \geq \alpha_0,
       \label{cond1}
    \end{equation}
    where $\pi:\hr\to\hh$ denotes the (vertical)
    projection onto the first factor.
  \end{claim}



    If not, for any sequence of points $p_n\in E$ such that
    $d(p_{n+1})>d(p_n)\geq n$ (the sequence diverges horizontally),
    the sequence $\{\nu(p_n)\}_{n\in \mathbb N}$ converges uniformly to zero. We then
    consider the isometry $T_n$ which is the composition of a vertical
    and a horizontal translations mapping $p_n$ to a fixed point
    $p_0=(0,0,\frac{L}{2})\in\h^2\times\R$, and call $E_n=T_n(E)$.

    If there is a subsequence of $\{p_n\}_{n\in \mathbb N}$ such that
    the Gaussian curvature of $E$ is uniformly bounded in a small
    neighborhood of any $p_n$, then there is a subsequence of the
    surfaces $E_n$ which locally converges to a CMC surface $E_\infty$ whose
    angle function $\nu_{\infty}$ vanishes identically in a neighborhood of $p$, and $E_\infty$
    is contained in a vertical cylinder over a complete curve of
    constant geodesic curvature~$2H$. By the unique continuation
    theorem, there is a subsequence of $\{E_n\}_{n\in \mathbb N}$
    converging to a complete vertical cylinder, a contradiction with
    the fact that all the terms in the sequence are contained in a
    fixed horizontal slab of height less than $3L$ ($\hh\times[-L, 2L]$).

    Then, we can assume that the curvature of $E$ is not uniformly
    bounded in neighborhoods of $p_n$ and, by passing to a
    subsequence if necessary, that $|A(p_n)|\geq n$ for any $n$, where $A$
    denotes the second fundamental form of $E$. We call $B_n$ the
    connected component that contains $p_n$ of the intersection of $E$
    with the extrinsic ball $\B(p_n,\delta)$ centered at $p_n$ of
    uniform small radius $\delta>0$, and we define
    \[
      f_n(p):= \overline{d}\left(p,\partial B_n\right) |A(p)|,
    \]
    for any
    $p\in B_n$, where $\overline{d}$ denotes the extrinsic distance in
    $\hr$.  The function $f_n$ vanishes on $\partial B_n$ and
    $f_n(p_n)=\delta|A(p_n)|\geq \delta n$.  We then deduce that $f_n$
    attains its maximum at a point $q_n\in B_n$, and
    $f_n(q_n)\geq \delta n$. On the other hand,
    $\delta|A(q_n)|\geq f_n(q_n)$, from where we deduce that
    $|A(q_n)|\geq n$.

    We now consider $r_n:=\frac 12 \overline{d}(q_n,\partial B_n)$ and
    $B'_n\subset B_n$ the connected component of $E\cap\B(q_n,r_n)$
    that contains $q_n$. For any point $q\in B'_n$, it holds
    \[
      2r_n=\overline{d}(q_n,\partial B_n)\leq \overline{d}(q_n,q) +
      \overline{d}(q, \partial B_n)\leq r_n +\overline{d}(q, \partial
      B_n).
    \]
    Hence, $\overline{d}(q,\partial B_n)\geq r_n$. Since
    $2r_n |A(q_n)|=f_n(q_n)\geq f_n(q)\geq r_n |A(q)|$, we conclude
    that
    \[
      |A(q)|\leq 2|A(q_n)|=:2\lambda_n.
    \]

    Since $\nu(p_n)$ converges uniformly to zero, we can assume, by
    passing to a subsequence if necessary, that $|\nu|<\frac 1n$ in
    $B_n$.  Thus we have that $|\nu|<\frac 1n$ and $|A|\leq 2\lambda_n$
    in $B'_n$.

    Now we consider a blow up on the metric $g$ of $\hr$ by a factor
    $\lambda_n\geq n$; more precisely, we define $\Sigma_n$ as $B'_n$
    with the metric $g_n=\lambda_n g$.  We can use the exponential map
    at the point $q_n$ to lift the surface $\Sigma_n$ to its tangent
    plane $T_{q_n}(\hr)\approx \R^3$, and we obtain a surface
    $\widetilde{\Sigma}_n\subset \R^3$ which is a minimal surface with
    respect to the lifted metric $\tilde{g}_n,$ where $\tilde{g}_n$ is
    the metric such that the exponential map exp$_{q_n}$ is an
    isometry from $(\widetilde{\Sigma}_n, \tilde{g}_n)$ to
    $(\Sigma_n,g_n)$. Therefore,
    $\widetilde{\Sigma}_n\subset \B_{\rr^3}(0,\lambda_{n}r_n)$ and, if
    $\widetilde A$ denotes the second fundamental form of
    $\widetilde{\Sigma}_n$, we have $|\widetilde A(0)|=1$ and
    $|\widetilde A(q)|\leq 2$ for all $q\in\widetilde{\Sigma}_n.$

    On one hand, we observe that $\lambda_n r_n$ diverges, as
    ${2\lambda_nr_n=f_n(q_n)\geq f_n(p_n)\geq\delta n}$.  Then the
    balls $\B_{\rr^3}(0,\lambda_nr_n)$ converge to $\R^3$ and the
    metrics $g_n$ converge to the canonical metric $g_0$ of $\R^3$.
    For a fixed $n$, the compact surfaces
    $\widetilde{\Sigma}'_k:=\widetilde{\Sigma}_k\cap
    B_{\rr^3}(0,\lambda_nr_n)$ of~$\R^3$, with $k\geq n$, all pass
    through the origin $0$ and have uniform bounded curvature. Then a
    subsequence of $\{\widetilde{\Sigma}'_k\}_k$ converges to a
    minimal surface in $(\rr^3,g_0)$ passing through the origin $0$
    and $|\widetilde A(0)|=1$, where $\widetilde A$ also denotes the
    second fundamental form of this limit surface. This argument holds
    for any $n$, so we can use a diagonal argument and obtain, as a
    limit of a subsequence of the surfaces $\widetilde\Sigma_k$, a
    complete minimal surface $\widetilde{\Sigma}$ in $\rr^3$ with
    $0\in \widetilde{\Sigma}$ and $|\widetilde A(0)|=1$.

    On the other hand, we knew that $|\nu|<\frac 1n$ in $B'_n$. Then
    we obtain that $|\nu|<\frac 1n$ in $\widetilde\Sigma_n$, from
    where we deduce that the Gauss map of $\widetilde\Sigma$ takes
    values in a neighborhood of the equator of $\s^2$. Then this limit
    surface $\widetilde\Sigma$ must be a vertical plane, which
    contradicts the fact that $|\widetilde A(0)|=1$. This proves
    Claim~\ref{cl}.

    \medskip
  
 
    Take the sequence $\{p_n\}_{n \in \N}$ given by Claim~\ref{cl}. We
    can assume, by passing to a subsequence if necessary, that
    $\nu (p_n)>0$ (or $\nu (p_n)<0$) for any $n$.  We consider, for
    each $p_n$, a translation of the paraboloid $\mathcal{P}^+$ (or
    $\mathcal{P}^-$, depending on the sign of $\nu (p_n) $) tangent
    to~$E$ at~$p_n$ with the same mean curvature vector at $p_n$. We
    denote such paraboloid by $\mathcal{P}(p_n)$.

    Since the angle function $\nu^+$ of the paraboloid $\mathcal{P}^+$
    is a decreasing function (by convexity of $u$, see Section
    \ref{sec-parab}), there exists a unique $\beta\in(0,\alpha_0)$
    such that $h_\beta=h_{\alpha_0} +L$. Hence,
    $\mathcal{P} (p_n) \cap (\h^2 \times [0,L])$ is the translation of
    a subdomain of
    $\mathcal{P}^\pm \cap (\h^2 \times [0,\pm h_\beta])$, where the
    sign $\pm$ depends on the sign of the angle function $\nu$ at
    $p_n$. For $n$ large enough, $d(p_n)>2 k_\beta$ and the tangent
    paraboloid $\mathcal{P} (p_n)$ does not intersect $\partial E$.

    The local intersection of $E$ and $\mathcal{P}(p_n)$ near $p_n$
    consists of $k$ curves, with $k\geq 2$, meeting at an equal angle
    at $p_n$. We denote by $\Gamma_n=E\cap \mathcal{P}(p_n)$ the
    intersection of the two surfaces. Since the intersection of
    $\mathcal{P}(p_n)$ with the slab containing $E$ is compact,
    $\Gamma_n$ is also compact.

    The paraboloid $\mathcal{P}(p_n)$ divides $\h^2\times\R$ into two
    components: a mean-convex one, $\mathcal{P}_{int}(p_n)$, and a non
    mean-convex one, $\mathcal{P}_{ext}(p_n)$.  Since the intersection
    of $\mathcal{ P}_{int}(p_n)$ with a horizontal slab is compact,
    any component of $E\setminus\mathcal{ P}(p_n)$ contained in
    $\mathcal{ P}_{int}(p_n)$ is necessarily compact.

    \begin{claim}\label{cl2}
      There is no compact component $\Sigma$ of
      $E\setminus\mathcal{ P}(p_n)$ contained in
      $\mathcal{P}_{ext}(p_n)$ with boundary
      $\partial \Sigma \subset \mathcal{P}(p_n)$.
    \end{claim}

    Suppose by contradiction this is not true. Then we can find a
    vertically translated copy $\widetilde{\mathcal{P}}(p_n)$ of
    $\mathcal{P}(p_n)$ tangent to $\Sigma$ at a point $\widetilde p$
    such that $\Sigma$ is contained in the mean-convex side of
    $\widetilde{\mathcal{P}}(p_n)$. By the maximum principle, we get
    that the mean curvature vector of $\Sigma$ at $\widetilde p$
    points to the non mean-convex side of
    $\widetilde{\mathcal{P}}(p_n)$.  Let $\sigma (\mathcal{P})$ be the
    symmetric copy of $\mathcal{P}(p_n)$ with respect to
    $\h^2\times\{0\}$. We can translate $\sigma(\mathcal{P})$ so that
    it is tangent to $\widetilde{\mathcal{P}}(p_n)$ (and hence to
    $\Sigma$) at $\widetilde p$. Then $\Sigma$ and
    $\sigma(\mathcal{P})$ share the same mean curvature vector at
    $\widetilde p$ and $\Sigma$ is contained in the non mean-convex
    side of $\sigma(\mathcal{P})$, a contradiction by the maximum
    principle, and Claim \ref{cl2} is proved.

    \medskip

  Locally at $p_n$, the set $E \setminus \Gamma_n$ has at least four
  components with at least two of them contained in the mean-convex
  component $\mathcal P_{int}(p_n)$. We call $\Sigma_1$ and $\Sigma_2$
  two of these components.

 \begin{claim}\label{cl3}
   $\Sigma_1$ and $\Sigma_2$ can be connected by an arc
   $\gamma \subset E \cap \mathcal{ P}_{int}(p_n)$.
  \end{claim}

  Suppose by contradiction this is not the case. Thus $\Gamma_n$
  bounds at least two distinct connected components $R_1$ and $R_2$ of
  $E \setminus \mathcal{ P}(p_n)$ which are contained in $\mathcal{
    P}_{int}(p_n)$ whose boundaries meet at~$p_n$. We consider
  vertical translations $T_s\mathcal{ P}(p_n)$ of the paraboloid,
  where $T_s(p)=p+s \partial_t$ for $s \geq 0$, which foliate
  $\mathcal{ P}_{int}(p_n)$. Since $E\cap \mathcal{ P}_{int}(p_n)$ is
  compact, there is a last leaf $T_{s_1}\mathcal{ P}(p_n)$ of the
  foliation that meets $R_1$. Then $T_{s_1}\mathcal{ P}(p_n)$ and
  $R_1$ are tangent at a point~$q_{1}$ and $R_1$ is below
  $T_{s_1}\mathcal{ P}(p_n)$. By the maximum principle, the mean
  curvature vector of $E$ and the mean curvature vector of the
  paraboloid $T_{s_1}\mathcal{ P}(p_n)$ are opposite at $q_1$. The
  component $R_1$ separates $ \mathcal{ P}_{int}(p_n)$ into two
  connected components and the mean-convex one $R_1^+$ is compact.
  Observe that $R_1$ and $R_2$ are local graphs near $p_n$. Since the
  mean curvature vector of $E$ at $p_n$ is pointing into $\mathcal{
    P}_{int}(p_n)$ and $E$ is embedded, the component $R_2$ is
  completely contained in $R_ 1^+$ (see Figure \ref{figure8}) and the
  compact component of $ \mathcal{P}_{int}(p_n)\setminus R_2$ is not
  mean-convex. We then reach a contradiction applying the maximum
  principle with the last leaf $T_{s_2}\mathcal{ P}(p_n)$ of the
  foliation that meets $R_2$. This proves Claim \ref{cl3}.
  
    \begin{figure}[h]
   \centering \includegraphics[height=2.7cm]{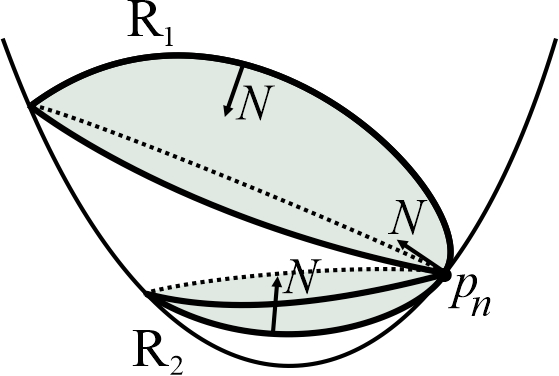}
   \caption{Connected components $R_1$ and $R_2$ of $E \setminus
     \mathcal{ P}(p_n)$ in $\mathcal{ P}_{int}(p_n)$ appearing in the
     proof of Claim~\ref{cl3}.}
  \label{figure8}
 \end{figure}

\medskip

Let $\gamma_n$ be a compact arc in $E \cap \mathcal{ P}_{int}(p_n)$  linking two points $q_1\in\Sigma_1$ and
$q_2\in\Sigma_2$. We can complete $\gamma_n$ by a compact segment
$\gamma'_i\subset\Sigma_i$ with endpoints $q_i$ and $p_n$ such that
$\alpha_n = \gamma_n \cup \gamma'_1\cup\gamma'_2$ is a loop in~$E$. If
$\alpha_n$ is homologous to zero in $E$, then it is the boundary of a
disk $D$ which contains points in $\mathcal{P}_{ext}(p_n)$ close to
$p_n$; hence, the disk $D$ has a subdomain in $\mathcal{P}_{ext}(p_n)$,
a contradiction to Claim~\ref{cl2}.  This proves that $\alpha_n$ is in
the homology class of $\partial E$ and $\alpha_n \cup \partial E$ bounds
a subannulus $\mathcal{A}_n$ of $E$.

Since we can do this construction for a sequence of diverging points
$p_n$, we can use a translation of the paraboloid $\mathcal{P}^+$ so
that $\mathcal{ P}^+ \cap \partial E=\emptyset$ and $\mathcal{A}_n$
and $\mathcal{P}^+$ are tangent at a point. Then by the maximum
principle we conclude that, along $\mathcal{A}_n$, the mean curvature
vector of $E$ points into the solid cylinder $\mathcal{A}_n^+$ bounded
by the annulus $\mathcal{A}_n$ and two disks with
boundaries~$\alpha_n$ and~$\partial E$.  We are then in the situation
of a mean convex cylinder with boundaries contained in two compact
subdomains $K(\alpha_n)$ and $K(\partial E)$ of $\hr.$ Since all the
curves $\alpha_n$ are contained in the intersection with the slab
${\cal S}_L$ with translated copies of the mean-convex side of the
same paraboloid, ${\cal P}^\pm\cap(\hh\times[0,\pm h_\beta])$, we can
suppose that the compact sets $K(\alpha_n)$ are compact balls of
uniform radius.

First assume $0<H<\frac 1 2$.  For $n$ large enough, we can suppose
that there exists a Scherk domain $\Omega$ bounded by arcs $A_1, A_2,
B_1, B_2$ with $\kappa(A_i)=2H$ and $\kappa(B_i)=-2H$, where $\kappa$
denotes the geodesic curvature with respect to the interior of
$\Omega$, such that $B_i$ separates $K(\alpha_n)$ from $K(\partial
E)$, for $i=1,2$ (see Figure \ref{figure22}). In particular, we can conclude that  the vertical 
planes $A_1\times\rr$ and $A_2\times\rr$ do not 
intersect the annulus $\mathcal{A}_n$; otherwise
 we could use hyperbolic translations of the CMC $H$-plane $A_i\times\rr$
 to get a contradiction with 
the maximum principle (notice that the mean curvature vector of 
$\mathcal{A}_n$ points into the compact region bounded by it, so we have 
the correct orientation to apply the maximum principle). Now moving
 the CMC $H$-planes $A_i\times\rr$ towards the annulus $\mathcal{A}_n$
 using hyperbolic translations, we can guarantee the existence of two
 disjoint vertical CMC $H$-planes $\Gamma_1\times\rr$ and
 $\Gamma_2\times\rr$ such that  the annulus
$\mathcal A_n$ is contained in the convex side of the two of them
 and  $\Gamma_i$ separates $A_1$ from $A_2$, for $i=1,2$.  We can then translate
downwards the complete Scherk $H$-graph over $\Omega$ so that it is
below the annulus ${\cal A}_n$. Now we get a contradiction using the
maximum principle and vertical translations of the complete Scherk
$H$-graph.


    \begin{figure}[h]
   \centering \includegraphics[height=4cm]{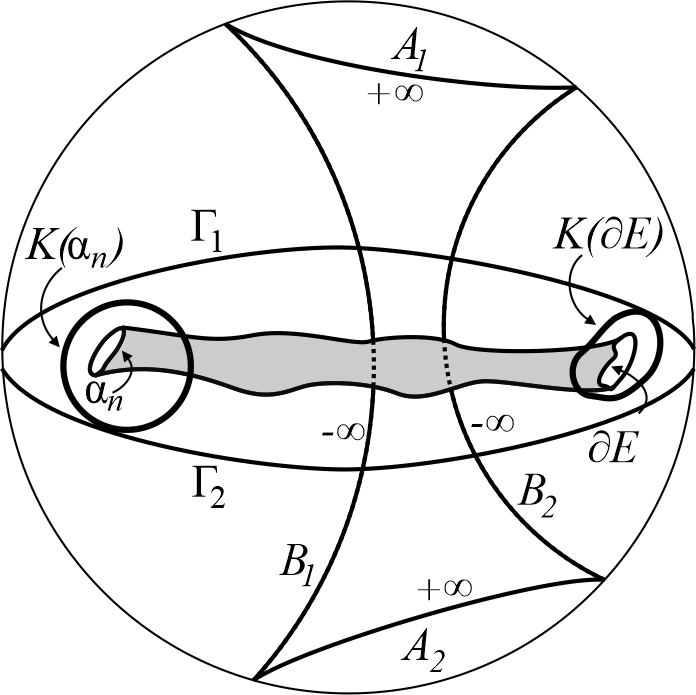}
\caption{Projection of the annulus ${\cal A}_n$ and of the complete
  Scherk $H$-graph for the case $0<H<\frac 12$.}
\label{figure22}

 \end{figure}

In the case of $H=\frac 1 2,$ the argument is similar but a bit more subtle and using
the graph given by the half of a horizontal catenoid described in
Section~\ref{sec.annuli} instead.

Let us first prove that, for any $n$, ${\cal A}_n$ projects onto the
region between two equidistant curves to a same geodesic.  Up to an
isometry, we can assume that the geodesic which minimizes the distance
between $K(\alpha_n)$ and $K(\partial E)$ is contained in
$\{y=0\}\subset\hh$, and let ${\cal I}$ be the entire $\frac 12$-graph
invariant by translations along this geodesic described in
Section~\ref{sec.annuli}. 

We consider $d>0$ so that the vertical projection of $K(\alpha_n)\cup
K(\partial E)$ is contained in ${\cal R}_d$, where we recall that
${\cal R}_d\subset\hh$ denotes the region bounded by the two
equidistant curves at distance $d$ to the horizontal geodesic
$\{y=0\}\subset\hh$.  We translate ${\cal I}$ downwards so that it is
much below ${\cal A}_n$ and start translating upwards until it
contains the equidistant curves $\partial{\cal R}_d\times\{L\}$ or
intersects $K(\alpha_n)\cup K(\partial E)$ for the first time. By the
maximum principle this translated copy of ${\cal I}$ lies below ${\cal
  A}_n$, and the vertical projection of ${\cal A}_n$ onto $\hh$ is
contained in ${\cal R}_d$. Since the compact balls $K(\alpha_n)$ have
uniform radius, this distance $d$ does not depend on $n$.

Let us now consider a horizontal catenoid $\mathcal C_a$ such that
${\cal R}_d$ is contained in the projection of ${\cal C}_a$ in $\hh$.
We can assume, up to a horizontal translation, that the vertical
symmetry plane $\Pi$ (which separates an end of ${\cal C}_a$ from the
other) intersects ${\cal A}_n$ transversally and has $\Pi\cap
K(\alpha_n)=\Pi\cap K(\partial E)=\emptyset$. For $M> L$ and $n$
sufficiently large, we can suppose that $\mathcal C_a\cap \{-M\leq
t\leq 0\}$ is a graph ${\cal G}$ over a domain $\tilde\Omega$ such
that $\tilde\Omega\cap K(\alpha_n)=\tilde\Omega\cap K(\partial
E)=\emptyset$ as in Figure~\ref{figure23}. Applying the maximum
principle with vertical translations of the graph ${\cal G}$, we
conclude that $\mathcal A_n$ can not exist.

    \begin{figure}[h]
   \centering \includegraphics[height=4cm]{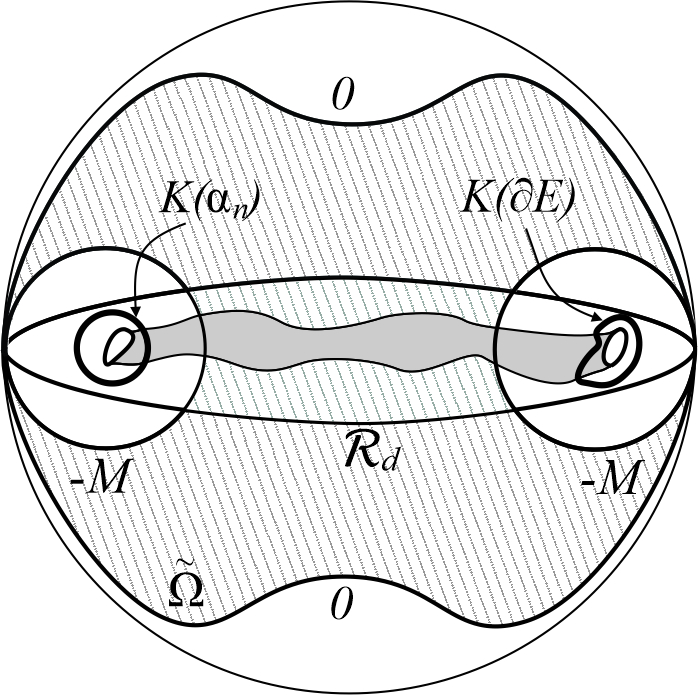}
\caption{Projection of the annulus ${\cal A}_n$ and of the graph
  ${\cal G}_a$ for the case $H=\frac 12$.}
\label{figure23}
 \end{figure}

\end{proof}

\section{Halfspace Theorem}

\begin{theorem}
  Let $M$ be a properly embedded CMC surface in
  $\hh\times [0,+\infty)$ with $0<H\leq \frac{1}{2}$ and finite
  topology. Then $M$ is necessarily a graph over a simply connected
  domain of $\hh$. For $H=\frac{1}{2}$ the graph is
  entire.
\label{thm2}
\end{theorem}

\begin{proof}
  Suppose by contradiction that there exist two distinct points $p_1,
  p_2\in M$ such that $\pi(p_1)=\pi(p_2)=x$, where $\pi:\hr\to\hh$ is
  the (vertical) projection onto the first component; that is,
  $p_1=(x,t_1)$ and $p_2=(x,t_2)$, with $t_1\neq t_2$. Assume
  $t_1<t_2$ and take $t_3>2t_2$ so that there exists a curve in $M$
  joining $p_1,p_2$ contained in $\{t<t_3\}$. We can assume that
  $M\cap\{t=t_3\}$ is transversal. Let us denote $\Sigma=M \cap
  \{t<t_3\}$.

  Suppose that $\overline \Sigma$ is compact. Then there exists a
  minimum for the height function~$t$, where the mean curvature vector
  $\overrightarrow{H}$ of $\Sigma$ coincides with~$\partial_t$. Hence,
  we can start the Alexandrov reflection method for $\Sigma$ with
  horizontal planes coming from below, and obtain that $M$ is
  symmetric with respect to a horizontal plane below $\{t=t_2\}$,
  which implies that $M$ is compact, a contradiction. Thus $\overline
  \Sigma$ cannot be compact.

  For any $\rho>0$, consider a translated paraboloid
  $T_t\mathcal{P}^-$ with $T_t (p) =p+t \partial_t$ so that the solid
  cylinder $\overline{D_{\hh}(x,\rho)}\times[0,t_3]$ is contained in
  $T_t\mathcal{P}^-_{int}$, where $D_{\hh}(x,\rho)$ is the hyperbolic
  disk centered at $x$ with radius $\rho$.  Since $M$ has finite
  topology, we can assume that, for $\rho >0$ large enough, any
  connected component of $M_\rho=M\cap T_t\mathcal{P}^-_{ext}$ is
  either compact with boundary on $T_t\mathcal{P}^-$ or an annular end
  of $M$ with boundary on $T_t\mathcal{P}^-$. Since $M$ is properly
  embedded, the number of compact components of $M_{\rho}$ is finite.

  For $\rho'>\rho$ large enough, we consider the vertical
  cylinder~$C_{\rho'}= D_{\hh}(x,\rho')\times \R$ which contains
  entirely $T_t\mathcal{P}^-\cap \mathcal{S}_{t_3}$ and all compact
  components of $M_\rho$ with boundary in the paraboloid, where we
  recall that $\mathcal{S}_{t_3}$ denotes the slab
  $\hh\times[0,t_3]$. If a point of $M$ is outside this cylinder, then
  it is part of an end of $M$ which has its boundary on the
  paraboloid~$T_t\mathcal{P}^-$. We will only consider the part of
  $M_\rho$ contained in the slab $\mathcal{S}_{t_3}$ and outside the
  cylinder~$C_{\rho'}$. Let us consider
  $R_{\rho'} =\mathcal{S}_{t_3}\setminus C_{\rho'}$.

  \begin{claim}\label{cl4}
    We can take $\rho'>\rho$ large enough so that
    $M_{\rho} \cap {\partial}C_{\rho'}\cap \mathcal{S}_{t_3} $ is transversal
    and any non-compact component of $M_{\rho}\cap R_{\rho'}$ satisfies $\nu >0$.
  \end{claim}

  Since $M$ is properly embedded, there are a finite number of
  non-compact components of $M_\rho$. It then suffices to prove the
  claim for any of them.  Let $E$ be a non-compact component of
  $M_\rho$ whose boundary is in $T_t\mathcal{P}^-$ and suppose by
  contradiction that there exists a diverging sequence of points
  $\{p_n\}_{n\in \mathbb N}$ in $E\cap R_{\rho'}$ with $\nu (p_n) \leq 0$.  If
  $\nu \equiv 0$ in a neighborhood of a point $p_n$, $M$ would be
  (locally at $p _n$) a cylinder over a curve of constant geodesic
  curvature $2H$. By the analytic continuation theorem, $M$ would be a
  complete CMC cylinder contradicting the fact that $M$ is contained
  in a halfspace. Then in a neighborhood of any $p_n$, there must
  exist a point $q_n\in M$ where $\nu(q_n)<0$. Hence, we can work
  with a sequence of points $q_n\in E\cap R_{\rho'}$ with
  $\nu (q_n) < 0$ and the sequence $\{q_n\}_{n\in\mathbb N}$ diverges
  to $(\partial_\infty\h^2)\times\R$.
  
  Now we are going to argue as in the proof of Theorem \ref{thm1}. We
  observe that the third coordinate of the points $q_n$ is bounded by
  $t_3$.  For any $n,$ since $\nu(q_n)<0,$ we can consider a
  translation of the paraboloid ${\mathcal P}^-$, denoted by
  ${\mathcal P}(q_n)$, tangent to $E$ at $q_n$ with same mean
  curvature vector.  A similar argument as in Claim~\ref{cl} yields
  that there exists a subsequence of divergent points
  $\{q_n\}_{n\in\mathbb N}$ such that $\nu(q_n)\leq \alpha_0<0,$ for
  some $\alpha_0\in \rr^{-}$, and the paraboloid ${\mathcal P}(q_n)$
  does not intersect $\partial E \subset T_t{\mathcal P}^-$ for $n$
  large enough. Now arguing as in Claims~\ref{cl2} and~\ref{cl3}, we
  can conclude that there is curve $\alpha$ in the homology class of
  $\partial E$ contained in the adherence of the mean-convex component
  determined by the tangent paraboloid. Hence, there is a subannulus
  $\cal A$ of $E$ with boundary $\alpha\cup\partial E$ and we can use
  a complete Scherk $H$-graph if $0<H<\frac 12$ or half a horizontal
  catenoid ${\cal C}_a$ if $H=\frac 12$, to reach a
  contradiction. This proves that for $\rho' >0$ large enough, any
  point $p$ of $M_\rho$ contained in $R_{\rho'}$ satisfies $\nu (p)
  >0$.  To conclude Claim \ref{cl4}, it suffices to take $\rho'$
  larger, if necessary, so that the intersection $M_{\rho} \cap
  \partial C_{\rho'}\cap S_{t_3}$ is transversal.

  \medskip

  We denote by $\widetilde \Sigma=M \cap C_{\rho'}\cap S_{t_3}$, with
  boundary
  $\partial \widetilde \Sigma \subset \partial C_{\rho'} \cup \{t=t_3\}$.
  We observe that if
  $ p= (y_0,t_0)\in \partial \widetilde \Sigma \cap \partial C_{\rho'}$
  then there is no point
  $(y_0,t)\in \partial \widetilde \Sigma \cap \partial C_{\rho'}$ with
  $t_0<t<t _3$, since $\nu (y_0,t_0)>0$ and any other point would
  also satisfy $\nu (y_0,t)>0$, contradicting the embeddedness property
  of the surface. Therefore,
  $\partial \widetilde \Sigma \cap \partial C_{\rho'} $ consists of curves
  projecting graphically onto $\partial D_{\hh}(x,\rho')$ and/or
  curves in $D_{\hh}(x,\rho')\times \{t_3\}$.

  Now we apply Alexandrov reflection method to $\widetilde \Sigma$.  We
  observe that the minimum for the height function on $\widetilde \Sigma$
  can be attained at an interior point and/or a boundary point
  in~$\partial C_{\rho'}$. On both cases (taking $\rho'$ larger if
  necessary) $\widetilde \Sigma$ is locally a graph near that points, and we
  can start the Alexandrov reflection for $\widetilde \Sigma$ using the
  family of horizontal planes coming from below. 
  We obtain that $\widetilde \Sigma$ must be symmetric with respect to
  a horizontal plane below $\{t=t_2\}$, but this is not possible as
  either $M$ would be compact (if $\partial\widetilde \Sigma\cap
  \partial C_{\rho'}$ is empty) or $\partial\widetilde \Sigma\cap
  \partial C_{\rho'}$ would contain a point with positive angle
  function $\nu$, contradicting Claim~\ref{cl4}.

  This contradiction proves that no two distinct points in $M$ have
  the same vertical projection. Standard arguments show that the mean
  curvature vector of $M$ is then nowhere horizontal, so $M$ is a
  complete multigraph. Manzano and the third author proved
  in~\cite{MR} that any complete multigraph with $0<H\leq \frac 1 2$
  (see \cite{HRS} for the case $H=\frac 1 2$) must be a graph over an
  unbounded domain $\Omega\subset\hh$.

  Now let us prove that $\Omega$ is simply connected. If this were not
  the case, there would be at least a Jordan curve $\gamma$ in
  $\partial\Omega$. Since $M$ is complete, the graph would necessarily
  diverge over $\gamma$. It was proved in~\cite{HRS2} that this is only
  possible if $\gamma$ has constant geodesic curvature $\pm 2H$, a
  contradiction since $\gamma$ is closed.

  Finally, observe that for $H=\frac 12$ the graph is necessarily
  entire by the result in \cite{HRS}.

\end{proof}

\Addresses

\end{document}